\documentclass[10pt,leqno]{article}
\usepackage{amssymb,amsbsy,amsmath,amsfonts,amssymb,amscd, mathrsfs}
\usepackage{color, graphicx} 

\usepackage[english]{babel}
\usepackage[T1]{fontenc}
\usepackage{indentfirst}

\makeatletter
\@addtoreset{equation}{section}
\makeatother

\newtheorem{statement}{}[section]
\newtheorem{theorem}[statement]{Theorem}
\newtheorem{lemma}[statement]{Lemma}

\newcommand\C{\mathbb C}
\newcommand\R{\mathbb R}

\newcommand\D{\mathbb D}

\newcommand\e{{\rm e}}
\newcommand\eps{\varepsilon}
\newcommand\ind{{\rm 1\kern-.30em I}}
\newcommand\qed{\hfill $\square$}

\let\phi=\varphi
\newcommand\capa{{\rm Cap}\,}
\newcommand\converge{\mathop{\longrightarrow}\limits}

\title{\bf A spectral radius type formula for approximation numbers of composition operators}
\author{\it Daniel Li, Herv\'e Queff\'elec, Luis Rodr{\'\i}guez-Piazza\footnote{Supported by a Spanish research project MTM 2012-05622.}}

\date{\footnotesize \today}

\begin{document}

\maketitle

\noindent{\bf Abstract.} \emph{For approximation numbers $a_n (C_\phi)$ of composition operators $C_\phi$ on weighted analytic Hilbert spaces, including 
the Hardy, Bergman and Dirichlet cases, with symbol $\phi$ of uniform norm $< 1$, we prove that 
$\lim_{n \to \infty} [a_n (C_\phi)]^{1/n} = \e^{- 1/ \capa [\phi (\D)]}$, where $\capa [\phi (\D)]$ is the Green capacity of $\phi (\D)$ in $\D$. This formula 
holds also for $H^p$ with $1 \leq p < \infty$.}

\medskip

\medskip


\section{Introduction} 

The determination of the approximation numbers of composition operators on Hilbert spaces of analytic functions on the unit disk (Hardy space, weighted Bergman 
space, Dirichlet space) is a difficult problem. Some partial results (see \cite{LIQUEROD}, \cite{LLQR}, \cite{estimates}, \cite{LELIQURO}, \cite{Herve-Seip}) 
show that no simple answer may be expected. However, we proved in \cite{LIQUEROD} and \cite{LELIQURO} that these approximation numbers cannot decay 
faster than geometrically: we always have $a_n (C_\phi) \geq c \, r^n$ for some constant $c > 0$ and some $0 < r < 1$. Moreover, we showed in those papers that 
$\lim_{n \to \infty} [a_n (C_\phi]^{1/n} = 1$ if and only if $\| \phi \|_\infty = 1$. \par
\smallskip

The quantity $\lim_{n \to \infty} [a_n (C_\phi)]^{1/n}$ looks like a spectral radius formula for the approximation numbers. Recall that if $T$ is a bounded 
operator on a complex Hilbert space $H$, with spectrum $\sigma (T)$, the classical spectral radius formula tells that for the spectral radius 
$r (T) := \sup_{\lambda\in \sigma(T)}|\lambda|$, one has the formula:
\begin{displaymath} 
r (T) = \lim_{n\to \infty} \Vert T^n\Vert^{1/n}
\end{displaymath} 
(the existence of the limit being part of the conclusion). \par

Now, if $a_n = a_{n} (T)$ is the $n$-th approximation number of a bounded operator $T$ on a Hilbert space $H$, it was shown (\cite{KON}, p.~133), by taking 
a rank-one perturbation of an $n$-dimensional shift, that, given $0 < \sigma < 1$, we can have $a_1 = \cdots = a_{n - 1} = 1$, and $a_n = \sigma$. Using 
orthogonal blocks of such normalized operators, one easily builds examples of compact operators $T$ for which the quantity  $[a_{n} (T)]^{1/n}$ has no limit 
as $n$ goes to infinity, and indeed satisfies:
\begin{displaymath} 
\liminf_{n\to \infty} [a_{n}(T)]^{1/n} = 0 \, , \qquad \limsup_{n\to \infty} [a_{n}(T)]^{1/n} = 1 \, .
\end{displaymath} 
We might as well use a diagonal operator with non-increasing positive diagonal entries $\eps_n$ such that $\liminf_n \eps_n^{\, 1/n} = 0$ and 
$\limsup_n \eps_n^{\, 1/n} = 1$. Nevertheless, the parameters
\begin{equation} 
\beta^{-} (T) = \liminf_{n\to \infty} [a_{n} (T)]^{1/n}, \qquad \beta^{+} (T) = \limsup_{n\to \infty} [a_{n}(T)]^{1/n}
\end{equation} 
which satisfy $0 \leq \beta^{-} (T) \leq \beta^{+} (T) \leq 1$ are similar to the term $\lim_{n \to \infty} \| T^n \|^{1/n}$ in the spectral radius formula. 
When the limit exists we will denote it by:
\begin{equation} 
\beta (T) = \lim_{n\to \infty} [a_{n} (T)]^{1/n} .
\end{equation} 
These parameters were shown to play an important role in the study of composition operators (see \cite{LIQUEROD} and \cite{LELIQURO}). As said above, 
the following was proved in these papers. 
\begin{theorem} \label{ancien} 
Let $H$ be a weighted Bergman space $\mathfrak{B}_\alpha$ (in particular the Hardy space $H^2$) or the Dirichlet space $\mathcal{D}$ and 
$\varphi \colon \D\to \D$ inducing a composition operator $C_\varphi \colon H \to H$. Then:
\smallskip

$1)$ if $0 < \Vert \varphi \Vert_\infty < 1$, one has $0 < \beta^{-} (C_\varphi) \leq \beta^{+} (C_\varphi) < 1$; \par
\smallskip

$2)$ if $\Vert \varphi \Vert_\infty = 1$, one has $\beta^{} (C_\varphi) = 1$. 
\end{theorem} 

The aim of this work is to complete this result by showing that $\beta (C_\varphi)$ exists as well when $\Vert \varphi\Vert_\infty < 1$ and to give a closed 
formula for this $\beta (C_\varphi)$  in terms of a Green capacity, relying on a basic work of \cite{WID} (see also \cite{FIMI}). We thus get another proof 
of $2)$ in the above theorem. \par
\smallskip

We end the paper with some words on the $H^p$ case for $1 \leq p < \infty$.
\par\medskip

We begin by giving notations, definitions and facts which will be used throughout this work. 

\section{Background, framework, and notations} 

Recall that if $X$ and $Y$ are two Banach spaces of analytic functions on the unit disk $\D$, and $\phi \colon \D \to \D$ is an analytic self-map of $\D$, one 
says that $\phi$ induces a \emph{composition operator} $C_\phi \colon X \to Y$ if $f \circ \phi \in Y$ for every $f \in X$; $\phi$ is then called the 
\emph{symbol} of the composition operator. One also says that $\phi$ is a symbol for $X$ and $Y$ if it induces a composition operator $C_\phi \colon X \to Y$.

\subsection{Singular numbers}

For an operator $T \colon X \to Y$ between Banach spaces $X$ and $Y$, its \emph{approximation numbers} are defined, for $n \geq 0$, as:
\begin{equation} \label{approx numbers} 
a_n (T) = \inf_{\text{rank}\, R < n} \| T - R\| \,.
\end{equation} 
One has $\| T \| = a_1 (T) \geq a_2 (T) \geq \cdots \geq a_n (T) \geq a_{n + 1} (T) \geq \cdots$, and (assuming that $Y$ has the Approximation Property), 
$T$ is compact if and only if $a_n (T) \converge_{n \to \infty} 0$.\par\smallskip

The \emph{$n$-th Kolmogorov number} $d_{n} (T)$ of $T$ is defined as (see \cite{CAST}, p.~49): 
\begin{equation} \label{Kolmo numbers} 
d_{n} (T) = \inf_{\substack{E \subseteq Y \\ \dim E < n}}\raise -2 pt \hbox{$\bigg[$} \sup_{x \in B_X} {\rm dist}\, (T x, E) \raise - 2 pt \hbox{$\bigg]$} 
= \inf_{\substack{E \subseteq Y \\ \dim E < n}} \| Q_E T \|_{Y/ E} \, , 
\end{equation} 
where $Q_E \colon Y \to Y / E$ is the quotient map. 
One always has $a_n (T) \geq d_n (T)$ and, when $X$ and $Y$ are Hilbert spaces, one has $a_n (T) = d_n (T)$ (see \cite{CAST}, p.~51). \par

\medskip
As usual, the notation $A \lesssim B$ means that there is a constant $c$ such that $A \leq C \, B$. 

\subsection{Weighted analytic Hilbert spaces} 

An \emph{analytic Hilbert space} $H$ on $\D$ is a Hilbert space $H \subset {\cal H}{\rm ol} (\D)$, the analytic functions on the unit disk $\D$, for which the 
evaluations $f \mapsto f (a)$ are continuous on $H$ for all $a \in \D$ and therefore given by a scalar product:
\begin{displaymath} 
f (a) = \langle f, K_a \rangle \, , \quad K_a \in H.
\end{displaymath} 
Since weakly convergent sequences of $H$ are norm-bounded, the \emph{reproducing kernels} $K_a$ are automatically norm-bounded on compact subsets of 
$\D$, that is:
\begin{equation}\label{automat} 
\qquad L_r := \sup_{|a|\leq r} \Vert K_a \Vert < \infty, \quad \text{for all } r < 1. 
\end{equation}

We will be slightly less general here, and adopt  the framework of \cite{KELE}. Let $\omega \colon [0, 1) \to (0,\infty)$ be a continuous, positive, and 
Lebesgue-integrable function. We extend this function to a radial weight on $\mathbb{D}$ by setting $\omega (z) = \omega (|z|)$. We denote by 
$H_\omega$ the space of analytic functions on $\mathbb{D}$ such that
\begin{displaymath} 
\Vert f \Vert_\omega^2 := |f (0)|^2 + \int_{\mathbb{D}}|f '(z)|^2 \,\omega (z) \, dA (z) <+ \infty,  
\end{displaymath} 
where $dA$ stands for the normalized area measure on $\D$. We will often omit the subscript $\omega$ and write $\Vert \, . \, \Vert$ for 
$\Vert \, . \, \Vert_\omega$.\par 

If $f (z) = \sum_{n = 0}^\infty b_n z^n$, a computation in polar coordinates shows that:
\begin{equation}\label{expression} 
\Vert f \Vert^2 = \sum_{n = 0}^\infty |b_n|^2 \, w_n \,
\end{equation}
where:
\begin{equation} 
\qquad w_0 = 1 \quad \text{and} \quad w_n = 2 n^2 \int_0^1 r^{2 n - 1} \omega (r) \, dr\, ,\quad  n\geq 1. 
\end{equation} 
Observe that there is a constant $C = C (\omega) \geq 1$ and, for each $\varepsilon > 0$, a $\delta_\varepsilon > 0$ such that:
\begin{equation}\label{siz}
\qquad \quad \delta_\varepsilon \, \e^{-\varepsilon n} \leq w_n \leq C \, n^2 \, , \qquad n\geq 1.
\end{equation}
Indeed, in one side, one has $w_n \leq 2 n^2 \int_0^1\omega (r) \, dr$, and, on the other side, for each $0 < \delta < 1$,  setting 
$c_\delta = \inf_{0 \leq r \leq \delta}\omega (r)$, we have $c_\delta > 0$ and:
\begin{displaymath} 
w_n \geq 2 n^2 \, c_\delta \int_{0}^\delta r^{2n  -1} \, dr = c_\delta\, n\,\delta^{2n} , 
\end{displaymath} 
giving \eqref{siz}. This shows in passing that $H_\omega$ is an analytic Hilbert space, and we call it a \emph{weighted analytic Hilbert space}. This 
framework is sufficiently general for our purposes and includes for example the case of (weighted) Bergman, Hardy, and Dirichlet spaces, corresponding to 
$\omega (r) = (1 - r^2)^\alpha$, $\alpha > - 1$, that is $w_n \approx n^{1 - \alpha}$. The standard  Bergman, Hardy, Dirichlet spaces correspond  to the 
respective values $\alpha = 2, 1, 0$.
\smallskip

The following simple fact will be used. Let  $a \in \mathbb{D}$ and $j\geq 0$; then:
\begin{equation}\label{deriv}  
f \mapsto f^{(j)} (a)  \text{ is a continuous linear form on } H. 
\end{equation} 
This holds for any analytic Hilbert space on $\mathbb{D}$, thanks to \eqref{automat}, and here can also be viewed as a consequence of \eqref{siz}.
\smallskip

An analytic self-map $\varphi \colon \mathbb{D}\to \mathbb{D}$ which induces a composition operator $C_\varphi \colon H \to H$ will be called a 
\emph{symbol} for $H = H_\omega$. In our space $H$, we have a quite easy case  for deciding if some $\varphi$ is a symbol.

\begin{lemma}\label{bol}  
If $\Vert \varphi \Vert_\infty < 1$, then $\varphi$ is a symbol \textnormal{if and only if} $\varphi \in H$. Equivalently, if and only if the \textnormal{positive} 
measure $\mu = |\varphi '|^2 \omega \, dA$ is finite. In that case, we moreover have 
$\| \phi^{k} \| \leq  C \, k\, \| \phi \|_\infty^k\, \| \phi \|$ for every $k \geq 1$.  
\end{lemma}

\noindent{\bf Proof.} If $\varphi$ is a symbol, then $\varphi = C_{\varphi} (z) \in H$. Conversely, let $\rho = \| \phi \|_\infty < 1$. We first note that, if 
$\phi \in H$, we have for any integer $k\geq 1$: 
\begin{equation}\label{puiss} 
\begin{split}
\Vert \varphi^k\Vert^2 
& = |\varphi (0)|^{2 k} + \int_{\mathbb{D}}\omega (z) \,  k^2 \, |\varphi (z)|^{2 (k - 1)} |\varphi ' (z)|^2 \, dA(z) \\
& \leq \rho^{2k} (1+ k^2 \rho^{- 2}) \, \Vert \varphi \Vert^2.
\end{split}
\end{equation}
Now,  let $\varepsilon > 0$ be such that $\rho \, \e^{\varepsilon} < 1$.  If $f (z) = \sum  b_{k} z^k \in B_H$, the unit ball of $H$, we have by 
\eqref{siz}: $|b_k| \leq w_{k}^{- 1/2} \leq C_\varepsilon \e^{k \varepsilon}$, so that, using \eqref{puiss}, we see that the series 
$\sum b_k \, \varphi^k = f \circ \varphi$ converges absolutely in $H$, which proves that $C_\phi$ is compact (and even nuclear). \qed

\subsection{Green capacity}

The \emph{Green function} $g \colon \mathbb{D} \times \mathbb{D} \to (0,\infty]$ of the unit disk $\D$ is defined as:
\begin{equation} 
g (z, w) = \log \Big\vert \frac{1 - \overline{w} z}{ z - w} \Big \vert \, .
\end{equation} 
If $\mu$ is a finite positive Borel measure on $\D$ with compact support in $\D$, its Green potential is: 
\begin{equation} 
G_\mu (z) = \int_\D g (z, w) \, d\mu (w) 
\end{equation} 
and its \emph{energy integral} is:
\begin{equation} 
I (\mu) = \iint_{\D \times \D}  g (z, w) \, d\mu (z) \, d\mu (w) \, .
\end{equation} 
Of course, 
\begin{equation} 
I (\mu) = \int_\D G_\mu (z) \, d\mu (z) \, .
\end{equation} 

For any subset $E$ of $\D$, one sets:
\begin{equation} \label{def V}
V (E) = \inf _\mu I (\mu)  \, , 
\end{equation} 
where the infimum is taken over all  probability measures $\mu$ supported by a compact subset of $E$. Then the \emph{Green capacity}\footnote{Actually the 
inner capacity, but for open and compact sets, it it is equal to the outer capacity and hence, \emph{is} the capacity: see \cite{Brelot}, Chapitre~V, p.~63. 
Choquet's Theorem (\cite{Choquet}; see also \cite{Brelot}, Chapitre~V, p.~66),  asserts that the inner capacity is equal to the outer capacity for all Borel sets.} 
of $E$ in $\D$ is:
\begin{equation} 
\capa (E) = 1 / V (E) \, .
\end{equation} 
If $K \subseteq \D$ is compact, the infimum in \eqref{def V} is attained for a probability measure $\mu_0$. If moreover $V (K) < \infty$ (i.e. $\capa (K) > 0$), 
this measure is \emph{unique} and is called the \emph{equilibrium measure} of $K$.  One always has $V (K) < \infty$ when $K$ has non-empty interior, since 
then $I (\lambda) < \infty$ where $\lambda$ is the normalized planar measure on some open disk $\Delta \subseteq K$. It is clear that we  have:
\begin{displaymath} 
K \subseteq L \Rightarrow V (K) \geq V (L) \Rightarrow \capa (K) \leq \capa (L) \, , 
\end{displaymath} 
i.e. $\capa (K)$ increases with $K$ and: 
\begin{displaymath} 
\capa (E) = \sup_{K \subseteq  E, K \text{ compact}} \capa (K) \, .
\end{displaymath} 

We refer to \cite{Brelot} and \cite{Conway} and to the clear presentation of \cite{NieSak} for the definition of the Green capacity and of its basic properties. 
Actually, in \cite{Brelot}, the capacity is defined by another way (see \cite{Brelot}, Chapitre~V, pp.~52--55), as follows.
\begin{lemma} \label{def capa equiv} 
For every compact set $K \subseteq \D$, one has:
\begin{align*}
\capa & (K)  \\
& = \sup \{ \| \mu\| \, ; \ \mu \text{ positive Borel measure supported by } K 
\text{ and } G_\mu \leq 1 \text{ on } \D\} 
\end{align*}
\end{lemma}
This is the definition of de la Vall\'ee-Poussin. Since our main result is based on H. Widom's paper \cite{WID}, it must be specified that he also used this definition 
in \cite{WID}. \par

Let us note, though we will not use that, that we also have:
\begin{align*}
\capa (K)  
& = \inf \{ \| \mu \| \, ; \ \mu \text{ positive Borel measure on } \D \text{ and } G_\mu \geq 1 \text{ on } K\} \\ 
& = \inf \{ \| \mu \| \, ; \ \mu \text{ positive Borel measure on } \D \text{ and } G_\mu \geq 1 \ q.e. \text{ on } K\} \, ,
\end{align*}
where {\it q.e.} means: out of a set of null capacity. The equivalence between these two definitions is shown in \cite{NieSak}, Lemma~4.1 (see also 
\cite{Brelot}, Chapitre~XI, p.~140 and pp.~144--145). 
\bigskip

An important fact for this paper is well-known to specialists on the (Green) capacity. This fact, kindly communicated to us with its proof by A. Ancona 
(\cite{ANC}),  is as follows.
\begin{theorem} \label{alano} 
For every \emph{connected} Borel subset $E$ of $\D$ whose closure $\overline{E}$ is contained in $\D$, one has:
\begin{equation} 
\capa (E) = \capa (\overline{E}) \, . 
\end{equation} 
\end{theorem} 

For sake of completeness, we provide details for the reader. We begin with a definition: a subset $E$ of $\D$ is said to be \emph{thin} 
(in French: ``\emph{effil\'e}'') at $u \in \overline E$ if there exists a function $s$ which is superharmonic in a neighbourhood of $u$ and such that 
\begin{displaymath} 
s (u) < \liminf_{\substack{v \to u \\ v \in E}} s (v) \, .
\end{displaymath} 
We denote by $\tilde E$ the union of $E$ and of points in $\overline E$ at which $E$ is \emph{not} thin (it is known that $\tilde E$ is the closure of $E$ for 
the fine topology: see \cite{Conway}, Proposition~21.13.10). Then:
\begin{lemma} \label{enough} 
If $E$ is a \emph{connected} Borel subset of $\D$ whose closure $\overline E$ is contained in $\D$, one  has:
\begin{displaymath} 
\tilde E = \overline{E} \, .
\end{displaymath} 
\end{lemma}

\noindent {\bf Proof.} Lemma~\ref{enough} is an immediate consequence of the following result (see \cite{Brelot}, Chapitre~VII, Corollaire, p.~89).
\begin{theorem} [Beurling-Brelot] \label{mousket} 
Let $E \subseteq \mathbb{D}$ and $u \in \overline E$. If $E$ is thin at $u$, there exist \emph{circles} with center $u$ and arbitrarily small radius $> 0$ 
which do not intersect $E$.
\end{theorem}

Indeed, taking the previous result for granted, suppose that $E$ is thin at $u \in \overline{E}$, $u \notin E$, and let $v_0 \in E$, with $|v_0 - u| = d > 0$. 
The function $\rho \colon E \to \R$ defined by $\rho (v) = |v - u |$ takes the value $d$ as well as arbitrarily small values since $u \in \overline{E}$. By the 
intermediate value theorem,  it takes every value in $(0, d]$, contradicting Theorem~\ref{mousket}. This contradiction shows that 
$\overline{E} \subseteq \tilde E$, thereby ending the proof of Theorem~\ref{alano}. \qed
\medskip

Now, 
\begin{lemma} 
One has: 
\begin{displaymath} 
\capa (E) = \capa (\tilde E) \, .
\end{displaymath} 
\end{lemma}

\noindent {\bf Proof.} We know (Cartan's Theorem) that $\capa (\tilde E \setminus E) = 0$ (see \cite{Conway}, Theorem~21.12.14, and Proposition 21.13.10, 
or see \cite{Brelot}, Chapitre~VII, p.~86 and Chapitre~V, p.~57, or \cite{Papa}, Proposition~8.2 and Proposition~8.3). Since the capacity of Borel sets is easily 
seen (see \cite{Brelot}, Chapitre~V, p.~62, or \cite{LAN}, Chap.~II, \S\,1, p.~145) to be a subadditive set function, one gets 
$\capa (E) \leq \capa (\tilde E) \leq \capa (E) + \capa (\tilde E \setminus E) = \capa (E)$.  \null \hfill \qed
\par\bigskip

Throughout this paper, for convenience, we sometimes use the notation:
\begin{equation} 
M (E) := \e^{- 1/ \capa (E)} = \e^{- V (E)}.
\end{equation} 
\goodbreak

\section{Main result} 

The goal of this paper is to prove the following result. 

\begin{theorem} \label{principal} 
Let $H$ be a  weighted analytic Hilbert space with norm $\Vert \, . \, \Vert$. Let $\varphi \colon \mathbb{D} \to \mathbb{D}$ be a symbol for $H$, with 
$\overline{\varphi (\mathbb{D})}\subseteq \mathbb{D}$. Then 
\begin{displaymath} 
\lim_{n \to \infty} [a_{n} (C_\varphi)]^{1/n} =: \beta (C_\varphi) 
\end{displaymath} 
exists and the value of this limit is:
\begin{equation} \label{valeur beta} 
\beta (C_\varphi) = \e^{- 1/\capa [\phi (\D)]}.
\end{equation} 
\end{theorem}

Note that, by Theorem~\ref{alano}, $\capa [\phi (\D)] = \capa [\overline{\phi (\D)}]$, so Theorem~\ref{principal} will follow immediately from 
Theorem~\ref{above} and Theorem~\ref{below} below. \par
\medskip

The proof is based on two results of H. Widom (\cite{WID}). Though those theorems are in the $H^\infty$ setting, we will be able to transfer them to our 
Hilbertian setting. Before giving this proof, we will check the result ``by hand'' with an explicit example. 

\subsection{A very special test case} 

Before going into the proof of Theorem~\ref{principal}, we are going to illustrate it in a simple situation. \par\smallskip

Let $\varphi$ be a symbol acting on  $H = H^2$ with $\Vert \varphi \Vert_\infty < 1$. We know from \cite{LIQUEROD} that $\beta^{+} (C_\varphi) < 1$, 
and  for very special $\varphi$'s we will show directly, without appealing to Widom's results, that \eqref{valeur beta} holds. 
\begin{theorem} \label{exemple} 
Let $\varphi(z) = \frac{a z + b}{c z + d}$ be a fractional linear  function mapping $\mathbb{D}$ into $\mathbb{D}$, i.e. :
\begin{displaymath} 
\qquad | a |^2 + | b |^2 + 2\, | \overline{a} b - \overline{c} d | \leq | c |^2 + | d |^2 \quad  \text{and} \quad | c | \leq | d | \, .
\end{displaymath} 
Then $\beta (C_\varphi) = \exp \Big[- \frac{1}{\capa (K)}\Big]$.
\end{theorem}

The example $\varphi (z) = z/(2 z + 1)$ shows that one cannot omit the  condition $\vert c \vert \leq \vert d \vert$. \par
Recall that the pseudo-hyperbolic distance on $\D$ is defined by:
\begin{equation} \label{dist pseudo-hyp}
\qquad \rho (z, w) = \Big| \frac{z - w}{1 - \overline{z} w} \Big| \, , \qquad z, w \in \D \, .
\end{equation} 
We denote by $\Delta (w, r) = \{ z \in \D \, ; \ \rho (z, w ) < r \}$ the open pseudo-hyperbolic disk of center $w$ and radius $r$. \par
We have the following two facts (\cite{NieSak}, p.~3173 for the first one).
\begin{lemma} 
Let $L = \overline{\Delta} (w, r)$ be a closed pseudo-hyperbolic disk of  pseudo-hyperbolic radius $r$. Then:
\begin{equation} \label{capart} 
\capa (L) = \frac{1}{\log (1/r)} \, \cdot
\end{equation}
\end{lemma}
\begin{lemma} \label{two} 
Let $u, v \colon \mathbb{D} \to \mathbb{D}$ be univalent analytic maps such that $u (\mathbb{D}) = v (\mathbb{D})$. Then, $u=  v \circ \psi$ where 
$\psi \in {\rm Aut}\, (\mathbb{D})$.
\end{lemma}

Indeed , by hypothesis $u = v \circ \psi$ with $\psi$ well-defined and holomorphic for $v$ is injective. Moreover, $u (\D) = v [\psi (\D)] = v (\D)$, whence 
$\psi (\D) = \D$, again because $v$ is injective. Finally $\psi$ is injective since $u$ is. \qed
\medskip

\noindent{\bf Proof of Theorem~\ref{exemple}.} We may assume $\| \phi \|_\infty < 1$. We first consider the case 
$\phi (z) = a z$, with $| a | < 1$. In that case, it is clear that $a_{n} (C_\varphi) = \vert a \vert^{n - 1}$, and hence $\beta (C_\varphi) = | a |$ 
and $\overline{\phi (\D)} = \overline{D} (0, | a |) = \overline{\Delta} (0, | a |)$. So that \eqref{valeur beta} holds in view of \eqref{capart}. \par\smallskip

In the general case, one might say that the conformal invariance of  $\capa$ and $\beta$ does the rest. Let us provide some details. \par

In general, $\phi (\D)$ is an euclidean disk, therefore a pseudo-hyperbolic disk $\Delta (w, r) := \{z \, ; \ \rho \, (z, w) < r \} = \psi_1 [\Delta (0, r)]$, where 
$\rho$ is the pseudo-hyperbolic distance and $\psi_1 \in {\rm Aut}\, (\D)$; one has the same radius since automorphisms preserve $\rho$.  If $h (z) = r z$, one 
therefore has $\phi (\D) = \psi_{1} [h (\D)]$ (since $\overline{\Delta} (0,  r)$ and the  euclidean disk $\overline{D} (0, r)$ coincide). From  Lemma~\ref{two}, 
$\phi = \psi_{1} \circ h \circ \psi_2$  with $\psi_2 \in {\rm  Aut}\, (\D)$, and so $C_\phi = C_{\psi_2} C_h C_{\psi_1}$, implying 
\begin{displaymath} 
\beta (C_\phi) = \beta (C_h) = r \, ,
\end{displaymath}  
by the ideal property. Moreover, 
\begin{displaymath} 
\capa [\phi (\D)] = \capa [h (\D)]
\end{displaymath}
by conformal invariance. Since we know that the desired equality between $\beta$ and $\capa$ holds  for $h$, we get the result. \qed 
\bigskip 

Let us \emph{numerically} test the claimed value of $\beta (C_\phi)$ on the affine example $\phi (z) = \phi_{a, b} (z) = a z + b$  with $a, b > 0$ and 
$a + b < 1$ (note that $C_{\phi_{a, b}}$ and $C_{\phi_{| a |, | b |}}$ are unitarily equivalent and have the same approximation numbers $a_n$, so that there is 
no loss of generality by assuming $a, b > 0$). In that case, the $a_n (C_\varphi) = a_n$ were computed \emph{exactly} by Clifford and Dabkowski 
(\cite{CLIDAB}). Their result is as follows. One sets:
\begin{equation} \label{clida-delta}
\Delta = (a^2 - b^2 - 1)^2 - 4 b^2  \quad \text{and} \quad Q = \frac{1 + a^2 - b^2 - \sqrt{\Delta}}{2 a^2} \, \cdot
\end{equation} 
Then, one has $a_n = a^{n - 1} Q^{n - 1/2}$, and so: 
\begin{equation} \label{clida}
\beta (C_\varphi) = a Q \, . 
\end{equation} 
\par

The  result of the  theorem can be tested on that  example. Indeed, we have $K := \overline{\varphi (\D)} = \overline{D} (b, a)$, so that 
(\cite{LAN}, p.~175--177):  
\begin{displaymath} 
\capa (K) = \frac{1}{\log \lambda} \, \raise 1 pt \hbox{,}
\end{displaymath} 
where $\lambda > 1$ is the biggest root of the quadratic polynomial 
\begin{displaymath} 
P (z) = a z^2 - (1 + a^2 - b^2) z + a \, . 
\end{displaymath} 
In explicit terms: 
\begin{displaymath} \label{koff} 
\e^{- 1/ \capa (K)} = \frac{1}{\lambda} = \frac{1 + a^2 - b^2 - \sqrt{\Delta_0}}{2 a} \, \raise 1pt \hbox{,} 
\end{displaymath}
with:
\begin{equation} \label{koff-delta} 
\Delta_0 = (1 + a^2 - b^2)^2 - 4 a^2 \, .
\end{equation} 
To get $\beta (C_\varphi) = \e^{- 1/ \capa (K)}$, it remains to compare \eqref{clida} and \eqref{koff}, using \eqref{clida-delta} and \eqref{koff-delta}, and 
to observe that 
\begin{displaymath} 
\Delta = \Delta_0 = (1 + a + b) (1 + a - b) (1 - a + b) (1 - a - b) \, .
\end{displaymath} 
%

\subsection{Widom's results reformulated}

We are going to state Widom's results in a form suitable for us. We first quote the following lemma from \cite{WID}.  
 
\begin{lemma} [Widom] \label{lemwid} 
Let  $K \subseteq \mathbb{D}$ be compact. Then, given $\varepsilon > 0$, there exists a cycle $\gamma$, which is a finite union of disjoint Jordan curves 
contained in $\D$, and whose interior $U$ contains $K$, and a rational function $R$ of degree $< n$, having no zero on $\gamma$ and all poles on 
$\partial \D$, such that, for $n$ large enough: \par

$1)$  $|R (z)| \geq \e^{- \varepsilon n}$ for $z \notin U$; \par\smallskip

$2)$ $ |R (z)| \leq \e^{\,\varepsilon n} \, \e^{- n/\capa (K)}$ for $z \in K$. 
\end{lemma}

The first theorem of Widom (\cite{WID}, Theorem~2,  p.~348), in which $\mathcal{C} (K)$ denotes the space of complex, continuous functions on $K$ with 
the sup-norm, can now be rephrased as follows. 
\begin{theorem} [Widom] \label{WID} 
Let $K\subseteq \D$ be a compact set, and $\eps > 0$. Then, there exist a constant $C_\eps > 0$ and, for every integer $n$ large enough, a rational function $R$ 
with poles on $\partial \D$ and points $\zeta_i \in \D \setminus K$  such that for every $g \in H^\infty$, one has:
\begin{equation}\label{harold} 
\Vert g - h\Vert_{\mathcal{C} (K)} \leq  C_\varepsilon \, \e^{\, \varepsilon n} \, \e^{- n/\capa (K)} \, \Vert g \Vert_\infty \, ,
\end{equation} 
where:
\begin{displaymath} 
\qquad h (w) = R (w) \sum_{\substack{i, k \\ 1 \leq k \leq m_i}} \hskip - 0,5 em c_{i, k} (g) \, (w - \zeta_i)^{- k} \quad \text{with} 
\quad  \sum_{i} m_i < n 
\end{displaymath} 
and the maps $g \in H^\infty \mapsto c_{i, k} (g)$ are linear. \par
Moreover, if $H$ is a weighted analytic Hilbert space, these maps, restricted to $H^\infty \cap H$, extend to continuous linear forms on $H$. 
\end{theorem} 

Widom's theorem precisely says the following. If $R$ and $\gamma$ are the  rational function and cycle of Lemma~\ref{lemwid},  let $\zeta_i$ be the zeros of $R$ 
inside $\gamma$. Consider, for $w \in K$, the function  
\begin{displaymath}
G (w)  = R (w) \, \bigg[\frac{1}{2 \pi i} \int_{\gamma} \frac{g (\zeta)}{R (\zeta) \, (\zeta - w)} \, d\zeta \bigg] \, ; 
\end{displaymath}
then, by the residues theorem, 
\begin{displaymath} 
G (w) = g (w) - R (w) \sum_{i, k} c_{i, k} (g) \, (w - \zeta_i)^{- k} = g (w) - h (w) \, , 
\end{displaymath} 
and Widom's theorem says that $\| G \|_{{\cal C} (K)} \leq  C_\varepsilon \, \e^{2\, \eps n} \, [M (K)]^n \, \Vert g \Vert_\infty$. \par
\smallskip
 
The only additional remark made here is that the $c_{i, k}$ are of the form 
\begin{displaymath} 
c_{i, k} (g) = \sum_{j \leq m_i - k} \lambda_{i, j, k} \, g^{(j)} (\zeta_i)
\end{displaymath} 
where $\lambda_{i, j, k}$ are fixed scalars, so that by \eqref{deriv} they extend to continuous linear forms on $H$. \par
\smallskip

Observe that the linear forms 
$g \mapsto c_{i, k} (g ')$ are also continuous on $H$ since 
\begin{equation}\label{useful} 
c_{i, k}(g ') = \sum_{j \leq m_i - k}\lambda_{i, j, k} \, g^{(j + 1)} (\zeta_i) \, . 
\end{equation}
This observation will be useful later.
\bigskip

\noindent {\bf Remark.} The rational function $h$ above is analytic in $\mathbb{D}$. Indeed, since the $\zeta_i$ are zeros of $R$, the polar factors 
$(w - \zeta_i)^{- k}$ are compensated by $R (w)$ with the right multiplicity, so that the only poles of $R$ have modulus $\geq 1$. However (see \cite{WID}, 
Lemma~1, p.~346), the poles of $R$ are located on $\partial{\mathbb{D}}$, but we cannot ensure that $h \in H$. Fortunately,  we will see that 
$h \circ \varphi \in H$, and this will be sufficient for our purposes. 
\bigskip

We will need  a second theorem of H. Widom (\cite{WID}, Theorem~7, p.~353), which goes as follows.
\begin{theorem} [Widom] \label{joe} 
Let $K$ be a compact subset of  $\mathbb{D}$ and $\mathcal{C} (K)$ be the space of continuous functions on $K$ with its natural norm. Set:
\begin{displaymath} 
\delta_{n} (K) = \inf_{E} \raise - 2pt \hbox{$\bigg[$} \sup_{f \in B_{H^\infty}} {\rm dist} \, (f, E) \raise - 2pt \hbox{$\bigg]$} \, , 
\end{displaymath} 
where $E$ runs over all $(n - 1)$-dimensional subspaces of $\mathcal{C} (K)$ and 
${\rm dist} \, (f, E)  = \inf_{h \in E} \Vert f - h \Vert_{\mathcal{C}(K)}$. Then 
\begin{equation} 
\delta_{n} (K) \geq \alpha \, \e^{- n/\capa (K)}  
\end{equation} 
for some positive constant $\alpha$.
\end{theorem}
\goodbreak

\subsubsection{The upper bound} \label{section upper} 

\begin{theorem} \label{above} 
Let $H$ be an analytic weighted Hilbert space with norm $\Vert \, . \, \Vert$. Let $\varphi \colon \mathbb{D} \to \mathbb{D}$ be a symbol for $H$, such that 
$\Vert \varphi \Vert_\infty = \rho < 1$. Then:
\begin{displaymath} 
\beta^{+}(C_\varphi) := \limsup_{n \to \infty} [a_{n} (C_\varphi)]^{1/n} \leq \e^{- 1/ \capa [\overline{\varphi(\mathbb{D})}]}.
\end{displaymath} 
\end{theorem}

\noindent{\bf Proof.} Fix $\varepsilon > 0$  such that  $\rho\, \e^{\varepsilon} < 1$. \par \smallskip

If $f (z) =\sum_{k = 0}^\infty b_k z^k \in H$, let $g (z) := S_{l} f (z) = \sum_{k = 0}^{l - 1} b_k z^k$, with $l = l (n)$ be an integer to be adjusted. \par

\begin{lemma} \label{partialsum}  
We have:
\begin{displaymath} 
\Vert f \circ \varphi - g \circ \varphi \Vert \leq K_\varepsilon \rho^{l} \e^{\varepsilon l} .
\end{displaymath} 
\end{lemma}
\noindent{\bf Proof.} For $f (z) = \sum_k b_k z^k$, we have: 
\begin{align*}
\Vert f \circ \varphi - g\circ \varphi \Vert 
& = \Big\Vert \sum_{k = l}^\infty b_k \varphi^k \Big\Vert 
\leq\sum_{k = l}^\infty |b_k|\,\Vert \varphi^k\Vert \\
& \leq \bigg(\sum_{k = l}^\infty |b_k|^2 w_k \bigg)^{1/2} \bigg(\sum_{k = l}^\infty \Vert \varphi^{k}\Vert^2 w_{k}^{-1}\bigg)^{1/2} 
\leq K_\varepsilon \rho^{l} \e^{\varepsilon l} \, , 
\end{align*}
by using Cauchy-Schwarz inequality, the fact that $\Vert f \Vert \leq 1$, the inequalities \eqref{siz}, and a geometric progression. \qed
\medskip

Also, remark that we have, by the Cauchy-Schwarz inequality:
\begin{align*} 
\Vert (S_{l}f)' \Vert_\infty 
&\leq \sum_{k = 0}^{l - 1} k \,|b_k| 
\leq \bigg(\sum_{k = 0}^{ l- 1} |b_k|^{2} w_{k} \bigg)^{1/2} \bigg(\sum_{k = 0}^{l - 1} k^2 w_{k}^{- 1}\bigg)^{1/2} \\
& \leq \Vert f \Vert \, \bigg(\sum_{k = 0}^{l - 1} k^2 w_{k}^{- 1}\bigg)^{1/2} .
\end{align*} 
Therefore, using \eqref{siz}, we see that the linear map $S'_l : H \to H^\infty$, defined by $S'_l (f) = (S_l f)'$, is continuous with a norm less than 
$ (\sum_{k = 0}^{l - 1} k^2 w_{k}^{- 1})^{1/2} \leq K_\varepsilon \, \e^{\varepsilon l} $. \par \medskip

We now use Theorem~\ref{WID}, with $K = \overline{\phi(\D)} \subseteq \D$ (and for $n - 1$ instead of $n$). Set, for $ n \geq 2$, large enough:
\begin{displaymath} 
\qquad h_{1} (w) = R (w) \sum_{\substack{i, k \\ 1 \leq k \leq m_i}} c_{i, k} (g') (w - \zeta_i)^{- k} \quad \text{with}\quad \sum_{i} m_i < n - 1 \, .
\end{displaymath} 
Recall that $h_1$ is analytic in $\D$. Remark that $h_1$ depends linearly on $f$ and the map $f \mapsto h_1$ has a rank $< n - 1$. We denote by 
$I_1 \in {\cal H}{\rm ol}\, (\D)$ the primitive of $h_1$ taking the value $g[\phi (0)]$ at $\phi (0)$:
\begin{displaymath} 
I_1 (z) = \int_{\phi (0)}^z h_1 (u) \, du + g [\phi (0)] \, .
\end{displaymath} 
Next, define an operator $A$ of rank $< n$ on $H$ (the continuity of $A$ being justified by \eqref{useful}) by the formula:
\begin{equation} \label{choice} 
A (f) = I_1 \circ \varphi \, .
\end{equation} 
Note that, even if $I_1\notin H$, we easily check on the integral representation of the norm that $I_1 \circ \varphi \in H$ since we assumed $\varphi \in H$, i.e. 
(see Lemma~\ref{bol}) that $\varphi$ is a symbol. \par \smallskip

Assuming for the rest of the proof that $\Vert f \Vert\leq 1$, we have the following lemma. 
\begin{lemma} \label{ext}
We have:
\begin{displaymath} 
\Vert g \circ \varphi - I_{1} \circ \varphi \Vert \leq K_\varepsilon \, \e^{\eps (n - 1)} \, \e^{\eps l} \, \e^{- (n - 1) /\capa (K)}  .
\end{displaymath} 
\end{lemma}
\noindent{\bf Proof.} Since $\varphi \in H$ and since $h_1 = I'_1$ approximates $g'$ uniformly on $K$ and 
$\Vert g' \Vert_\infty = \| (S_l f)' \|_\infty \leq K_{\varepsilon} \e^{\varepsilon l}$, we have, by Theorem~\ref{WID}: 
\begin{align*}
\Vert g \circ \varphi - I_{1} \circ \varphi \Vert^2 
& =\int_{\mathbb{D}} \big|g' [\varphi (z)] - h_{1} [\varphi (z)] \big|^2 |\varphi ' (z)|^2 \omega (z) \, dA (z) \\
& \leq K_{\varepsilon}^{2} \e^{2 \eps (n - 1)} [M (K)]^{2 (n - 1)} \Vert g' \Vert_{\infty}^{2} \int_{\mathbb{D}}|\varphi ' (z)|^2 \omega (z) \, dA (z) \\
& \leq C \, K_{\varepsilon}^3 \e^{2\varepsilon l} \, \e^{2\eps (n - 1)} [M (K)]^{2 (n - 1)}  \, , 
\end{align*}
(with $C = \| \phi \|_\omega^2$), hence the lemma, provided that we increase $K_\varepsilon$. \qed 
\bigskip

We can now end the proof of Theorem~\ref{above}. \par\smallskip

Writing:
\begin{align*} 
\Vert C_{\varphi} (f) - A (f) \Vert 
& = \Vert f \circ \varphi - I_1\circ \varphi \Vert \\ 
& \leq \Vert f \circ \varphi - g\circ \varphi \Vert + \Vert g \circ \varphi - I_1 \circ \varphi \Vert \,, 
\end{align*} 
we have: \par\smallskip

1)  $\Vert f \circ\varphi - g \circ \varphi \Vert \leq K_\varepsilon \, \rho^l \, \e^{\varepsilon l}$ by Lemma~\ref{partialsum}; \par 
\smallskip

2) $\Vert g \circ \varphi - I_1 \circ \varphi \Vert \leq K_\varepsilon \e^{\eps (n - 1)} [M (K)]^{n - 1} \e^{\varepsilon l}$ by Lemma~\ref{ext}. \par
\medskip

\noindent We therefore get, since $a_n := a_{n} (C_\varphi) \leq \Vert C_\varphi - A \Vert$: 
\begin{displaymath} 
a_n \leq K_\varepsilon \, \rho^l \, \e^{\eps l} +  K_\eps \, \e^{\eps l} \e^{\eps (n - 1)} [M (K)]^{n - 1} .
\end{displaymath} 
Next, since $(a + b)^{1/n} \leq a^{1/n} + b^{1/n}$, we infer that:
\begin{equation}\label{laststep}
a_n^{1/n} \leq (K_\varepsilon)^{1/n} (\rho \, \e^\eps)^{l /n}  + K_{\eps}^{1/n} \e^{\eps l/n} \e^{\eps (n - 1)/ n} M (K)^{(n - 1)/ n}.
\end{equation}
We now adjust $l = N n$, where $N$ is a fixed positive integer, and pass to the upper limit with respect to $n$ in \eqref{laststep}. We get:  
\begin{displaymath} 
L := \limsup a_n^{1/n} \leq [\rho \, \e^{\varepsilon}]^{N}+ \e^{\varepsilon} \e^{\varepsilon N} M (K).  
\end{displaymath} 
Letting $\varepsilon$ go to $0$, we get $L\leq \rho^{N}+ M (K)$. Finally, letting $N$ tend to infinity, we get $L \leq M (K)$ as claimed, and that ends the proof of 
Theorem~\ref{above}. \qed

\subsection{The lower bound} 
\begin{theorem} \label{below}
Let $H$ be a weighted analytic Hilbert space and $\phi \in H$ such that $\| \phi \|_\infty < 1$. Then:
\begin{displaymath} 
\beta^{-} (C_\phi) := \liminf_{n \to \infty} [a_{n} (C_\phi)]^{1/n} \geq \e^{- 1 /  \capa [\phi (\D)] } \, .
\end{displaymath} 
\end{theorem} 

It will be convenient to work with the Kolmogorov numbers $d_n (C_\phi)$ instead of the approximation numbers $a_n (C_\phi)$. Recall that, for Hilbert spaces, 
one has $d_n (C_\phi) = a_n (C_\phi)$. We begin with a simple lemma, undoubtedly well known to experts, on approximation numbers of an operator $T$ on a 
Hilbert space $H$.  
\begin{lemma} \label{known} 
For every Hilbert space $H$ and every compact operator $T \! \colon \! H \! \to \! H$, one has, $B_H$ denoting the unit ball of $H$:
\begin{equation} \label{andre} 
d_{n} (T) = \inf_{\dim E < n} \raise -2pt \hbox{$\bigg[$} \sup_{f \in B_H} {\rm dist}\, \big(T f,  T (E) \big) \raise -2pt \hbox{$\bigg]$} \, .
\end{equation}
\end{lemma}
\noindent{\bf Proof.} Indeed, if $\varepsilon_{n}(T)$ denotes the right hand side in \eqref{andre}, we clearly have $d_n (T) \leq \varepsilon_n (T)$. Now, let:
\begin{displaymath} 
T f = \sum_{j = 1}^\infty a_j (T) \, \langle f, v_j\rangle \, u_j, 
\end{displaymath} 
with $(u_j)$ and $(v_j)$ two orthonormal sequences, be the Schmidt decomposition of $T$. Let $E_0$ be the span of $v_1,\ldots, v_{n - 1}$. Observe that 
$u_j = T (a_j^{- 1} v_j) \in T (E_0)$ for $j < n$. Now, if $f \in B_H$, one has:
\begin{align*} 
\big[{\rm dist}\, \big(T f, T( E_0) \big) \big]^2 
& = \bigg\Vert \sum_{j = n}^\infty a_j (T) \, \langle f, v_j \rangle \, u_j \bigg\Vert^2 
= \sum_{j = n}^\infty [a_{j} (T)]^{2}| \, \langle f, v_j \rangle|^2 \\
& \leq [a_{n} (T)]^{2} \sum_{j = n}^\infty | \langle f, v_j\rangle|^2 \leq [a_{n} (T)]^{2} ;
\end{align*} 
so that $\varepsilon_n (T) \leq \sup_{f \in B_H} {\rm dist}\, \big(T f, T (E_0) \big) \leq a_n (T) = d_n (T)$. \qed
\bigskip

\noindent {\bf Proof of Theorem~\ref{below}.} Let $0 < r_j < 1,$  $r_j \to 1$ and $\psi_j \colon \mathbb{D} \to \mathbb{D}$ be given by 
$\psi_{j} (z) = r_j \, z$. Set $K_j = \overline{\varphi \circ \psi_j (\mathbb{D})} = \overline{\varphi (r_j \mathbb{D})}$. Let $E$ be a 
subspace of $H$ of dimension $< n$. By restriction, $E$ can be viewed as a subspace of $\mathcal{C} (K_j)$. By the second result of Widom 
(Theorem~\ref{joe}), we can find $f \in B_{H^\infty}$, $f (z) = \sum_{k \geq 0} b_k z^k$, such that:
\begin{displaymath} 
\qquad \quad \Vert f - h \Vert_{\mathcal{C} (K_j)} \geq 2 \alpha \, [M (K_j)]^n \, , \qquad \forall h \in E \, ,  
\end{displaymath} 
where $\alpha > 0$ is an absolute constant. If $H^\infty$  contractively embeds into $H$, we can continue with this $f$. In the general case, we have to correct $f$ 
in order to be in $B_H$, the unit ball of $H$. To that effect, we simply consider a partial sum:
\begin{displaymath} 
g (z) = \sum_{k = 0}^{l - 1} b_k z^k  
\end{displaymath} 
and we note that, setting $\rho_j = \sup_{w \in K_j} |w|$, one has $\rho_j < 1$ and: 
\begin{align} 
\Vert f - g \Vert_{\mathcal{C} (K_j)} & \leq \frac{\rho_j^{\, l}\ }{(1 - \rho_j^2)^{1/2}} \label{UN}  \\
\Vert g \Vert_H & \leq C \, l \, , \label{DEUX} 
\end{align}
where $C = C (\omega) \geq 1$ is the constant appearing in \eqref{siz}. \par

Indeed, we have $\Vert f - g \Vert_{\mathcal{C} (K_j)} \leq \sum_{k = l}^\infty |b_k| \, \rho_j^k$ and then \eqref{UN} follows from Cauchy-Schwarz's 
inequality and the fact that $\sum_{k \geq 0} |b_k|^2 \leq 1$ since $f \in B_{H^\infty}$. For \eqref{DEUX}, we simply use that, by \eqref{siz}, the weight 
$w$ satisfies $w_k \leq C\, (k + 1)^2$ and get:
\begin{displaymath} 
\Vert g \Vert_H^2 = \sum_{k = 0}^{l - 1} \vert b_k \vert^2 \, w_k \leq  C \, l^2 \sum_{k = 0}^{l - 1} \vert b_k \vert^2 \leq C \, l^2 \leq C^2 \, l^2 .
\end{displaymath} 

We then notice that \eqref{UN} gives, for every $h \in E$: 
\begin{equation} \label{repair}
\begin{split}
\Vert g - h\Vert_{\mathcal{C} (K_j)}
\geq \Vert f & - h \Vert_{\mathcal{C} (K_j)} - \Vert f - g \Vert_{\mathcal{C} (K_j)} \\
& \geq 2 \, \alpha \, [M (K_j)]^n - \frac{\rho_j^{\, l} \ }{(1 -\rho_j^2)^{1/2}} \geq \alpha \, [M(K_j)]^n \, ,
\end{split}
\end{equation} 
if we take $l = A_{j}n$ where $A_j$ is a large positive integer depending only on $j$. Explicitly: 
\begin{displaymath} 
A_j > \frac{\log  \big[1 / \big(\alpha \, (1 - \rho_j^2)^{1/2} \big) \big]}{\log (1 / \rho_j)} + \frac{\log [ 1 / M (K_j) ] }{\log (1/\rho_j)} 
\, \cdot
\end{displaymath} 
\smallskip

Finally,  set $F = g / C \, l$. Then $F \in B_H$. Since $E$ is a vector space, \eqref{DEUX} and \eqref{repair} imply:
\begin{displaymath} 
\Vert F - h \Vert_{\mathcal{C} (K_j)} = \frac{1}{C \, l} \Vert g - C \, l \, h \Vert_{\mathcal{C} (K_j)} \geq \frac{1}{C \, l} \, \alpha \, [M (K_j)]^n .
\end{displaymath} 
But we also know that: 
\begin{displaymath} 
\Vert F - h \Vert_{\mathcal{C} (K_j)} 
= \Vert F \circ \varphi \circ \psi_j - h \circ \varphi \circ \psi_j \Vert_{\infty} 
\leq L_{r_j} \Vert F \circ \varphi - h \circ \varphi \Vert_H \, ,
\end{displaymath} 
so we are left with (recall that $l = A_{j}n$):
\begin{displaymath} 
\qquad \Vert C_\phi F - C_\phi h \Vert_H \geq \frac{\alpha}{C \, L_{r_j} A_j}\,\frac{M (K_j)^n}{n} \, \raise 1pt \hbox{,} \quad  \forall h \in E, 
\end{displaymath} 
implying by Lemma~\ref{known}:
\begin{displaymath} 
a_n (C_\varphi) = d_n (C_\varphi) \geq \frac{\alpha}{C \, L_{r_j} A_j} \,\frac{[M (K_j)]^n}{n} \, \cdot
\end{displaymath} 

Now, taking $n$-th roots and passing to the lower limit, we get:
\begin{equation} \label{presque} 
\beta^{-}(C_\varphi)\geq M(K_j).
\end{equation}

It remains now to let $j \to \infty$. Observe that the compact subsets $K_j \subseteq \phi(\D)$ form an exhaustive sequence of compact subsets of $\phi (\D)$. 
Let then $L \subseteq \phi (\D)$ be compact; we have $L \subseteq K_{j_0}$ for some $j_0$, and using \eqref{presque}, we get  
$\beta^{-} (C_\phi) \geq M (K_{j_0}) \geq M (L)$. Passing to the supremum on $L$, we get $ \beta^{-}(C_\phi) \geq M [\phi (\D)]$, and 
this ends the proof of Theorem~\ref{below}. \qed

\goodbreak

\subsection{The case $\| \phi \|_\infty = 1$} \label{section norme 1} 

As said in the Introduction, for weighted Bergman spaces (including the Hardy space), and for the Dirichlet space, we proved in \cite{LIQUEROD} and 
\cite{LELIQURO}, respectively, that $\beta (C_\phi) = 1$ if $\| \phi \|_\infty = 1$ for every $\phi$ inducing a composition operator on one of those spaces. 
\par
In this section, we use Theorem~\ref{principal} to generalize this result to all composition operators $C_\phi$ on weighted analytic Hilbert spaces, with another, 
and simpler, proof. \par

For that, it suffices to use the following result, which is certainly well-known to specialists. The pseudo-hyperbolic metric $\rho$ on $\D$ is defined in 
\eqref{dist pseudo-hyp} and we denote by ${\rm diam}_\rho$ the diameter for this metric. 
\begin{theorem} \label{saksman} 
Let $K$ be a compact and connected subset of $\D$. Then, for $0 < \varepsilon < 1$: 
\begin{displaymath} 
{\rm diam}_\rho\, K > 1 - \varepsilon \quad \Longrightarrow  \quad \capa (K) \geq c \log 1/\varepsilon \,, 
\end{displaymath} 
for some absolute positive constant $c$. \par 
Hence, the Green capacity of $K$ tends to $\infty$ as its pseudo-hyperbolic diameter tends to $1$.
\end{theorem} 
\medskip

Before proving that, let us give two suggestive examples, borrowed from  \cite{LAN}, p.~175--177. \par\smallskip

1) Let $K = \overline{D} (0, r)$; then:
\begin{displaymath} 
{\rm diam}_\rho \, K = \frac{2 r}{1 + r^2} \quad \text{and} \quad \capa (K) = \frac{1}{\log 1/r} \, \cdot
\end{displaymath} 
One sees that $r$ goes to $1$ when ${\rm diam}_\rho \, K$ goes to $1$, and hence $\capa (K)$ tends to infinity. 
\par\smallskip

2) Let $K = [0, h]$, with $0 < h < 1$. Then:
\begin{displaymath} 
{\rm diam}_\rho \, K = h \quad \text{and} \quad \capa (K) = \frac{1}{\pi}\, \frac{I'}{I} \, \raise 1pt \hbox{,} 
\end{displaymath} 
where $I$ and $I'$ are the elliptic integrals:
\begin{displaymath} 
I = \int_0^1 \frac{1}{\sqrt{(1 - t^2) (1 - k^2 \, t^2)}} \, dt \quad \text{and} \quad  I' = \int_0^1 \frac{1}{\sqrt{(1 - t^2)(1 - k'^{2}\, t^2)}} \, dt  \,, 
\end{displaymath} 
with $k = \frac{1 - h}{1 + h}$ and ${k'}^2 = 1 - k^2$. \par

If $0 \leq a < b \leq h$, then $b - a + h a b \leq h - a + a h^2 = h - a (1 - h^2) \leq h$, so that $\rho (a, b)\leq h$. Therefore, in this example again, the assumption 
${\rm diam}_\rho \, K \longrightarrow 1$ implies successively that $h \to 1$, $k \to 0$, $k' \to 1$, $I \to \pi/2$, $I' \to \infty$, and at last 
$\capa (K) \to \infty$. \par

This example shows that Theorem~\ref{saksman} is optimal since
\begin{displaymath} 
\int_0^1 \frac{dt}{\sqrt{(1 - t^2)(1 - k'^2 t^2)}} \approx \log \frac{1}{1 - {k'}^2} \approx \log \frac{1}{1 - h} 
\end{displaymath} 
as $h$ (and hence $k'$) goes to $1$.
\bigskip\goodbreak

The following proof of Theorem~\ref{saksman} was kindly shown to the second-named author by E. Saksman (\cite{EER}). \par 
It make use of the following alternative definition of Green capacity, where ${\cal C}_0^\infty (\D)$ is the space of infinitely differentiable functions on $\D$ 
which are null on $\partial \D$, and $d z = dx dy$ is the usual $2$-dimensional Lebesgue measure.
\begin{lemma} \label{altern} 
For every compact subset $K$ of $\D$, one has:
\begin{displaymath} 
\capa (K) = \inf\Big\{\frac{1}{2 \pi} \int_\D |\nabla u (z) |^2 \, dz \, ; \  u \in {\cal C}_0^\infty (\D) \text{ and } u \geq 1 \text{ on } K \Big\} \, .
\end{displaymath}
\end{lemma}

\noindent{\bf Proof  of Theorem~\ref{saksman}.} If ${\rm diam}_\rho \, K > 1 - \varepsilon$ and  $K$ is connected, it contains two points $z_1$ 
and $z_2$ such that $\rho (z_1, z_2) = 1 - \varepsilon$. By the invariance of the Green capacity and of $\rho$ under automorphisms of the disk, we can assume 
that $z_1 = 0$ and $z_2 = 1 - \varepsilon$. Take $\varepsilon < r < 1$.  Denote by $\Delta_r$ the intersection of the closed disk with center $1$ and radius $r$ 
with the closed unit disk. We observe that $K$ meets the exterior of $\Delta_r$ at $0$ and its interior at $1 - \varepsilon$. The connectedness of $K$ implies that 
$K$ meets the boundary of $\Delta_r$: there is $b \in K$ such that $|b - 1| = r$. Write $b = 1 + r \e^{i\vartheta}$. Take now $a = 1 + r \e^{i \theta}$ with 
$|a| = 1$ and $0 \leq \theta \leq \vartheta \leq 2 \pi$. Since $u (a) = 0$ and $u (b) \geq 1$, we get, by the fundamental theorem of calculus, that:
\begin{align*} 
1\leq u (b) - u (a) 
& = \int_\theta^\vartheta i r \, \e^{it} \nabla u (1 + r \e^{i t}) \, dt = \bigg| \int_\theta^\vartheta i r \, \e^{it} \nabla u (1 + r \e^{i t}) \, dt \bigg| \\ 
& \leq r \int_{\theta}^{\vartheta} |\nabla u (1 + r \e^{i t})| \, dt 
\leq  r \int_{0}^{2\pi} | \nabla u (1 + r \e^{i t})|\, dt.  
\end{align*} 
Now, Cauchy-Schwarz inequality gives: 
\begin{displaymath} 
\int_{0}^{2 \pi} |\nabla u (1 + r \e^{i t})|^2 \, dt \geq \frac{1}{2 \pi r^2} \, \cdot
\end{displaymath} 
Integrating in polar coordinates centered at $1$ and remembering that $u = 0$ outside $\mathbb{D}$, we get:
\begin{align*} 
\int_{\mathbb{D}} |\nabla u (z)|^2 \, dz 
& \geq \int_{\varepsilon < |z - 1 | < 1} |\nabla u (z)|^2 \, dz \\
& = \int_{\varepsilon}^{1} \bigg[\int_{0}^{2\pi}|\nabla u (1 + r \e^{i t})|^2 \, dt \bigg] \, r \, dr \\ 
& \geq \frac{1}{2 \pi} \int_{\varepsilon}^1 \frac{dr}{r} = \frac{1}{2 \pi} \log \frac{1}{\varepsilon} \, \cdot
\end{align*} 
In view of \eqref{altern}, this ends the proof of Theorem~\ref{saksman}. \qed
\par\bigskip\goodbreak

\noindent{\bf Proof of Lemma~\ref{altern}.} Though this result is often considered as ``well-known'', we were not able to find anywhere an explicit reference. 
Since the average reader (if any!) of this paper will not be a specialist in Potential theory, we give a proof here. \par \smallskip

1) We first prove that the capacity of the compact $K$ is less than the right-hand side (though we only need that it is greater). We shall use  
Lemma~\ref{def capa equiv}. \par

We know (\cite{Conway}, Corollary~21.4.7, or \cite{LAN}, p.~91) that for every measure $\mu$ on $\D$ supported by $K$, one has 
$\Delta G_\mu = - 2 \pi \mu$, where $G_\mu$ is seen as a distribution. Hence, for every function $u \in {\cal C}_0^\infty (\D)$ such that $u \geq 1$ on $K$ 
and every positive measure $\mu$ supported by $K$ such that $G_\mu \leq 1$ on $\D$, one has:
\begin{displaymath} 
\mu (K) = \int_K d \mu \leq \int_\D u \, d \mu = - \frac{1}{2 \pi} \int_\D u (z) \, \Delta G_\mu (z) \, dz \, .
\end{displaymath} 
Then, by definition of the Laplacian of a distribution, we get:
\begin{displaymath} 
\mu (K) \leq - \frac{1}{2 \pi} \int_\D \Delta u (z) \, G_\mu (z) \, dz \, . 
\end{displaymath} 
But (see \cite{Brelot}, Chapitre~XI, p.~132 and pp.~144--145, or \cite{LAN}, Chap.~IV, \S~1, p.~215), for every real Borel measures $\nu_1$ and $\nu_2$ 
with finite energy (meaning that their positive and negative parts have finite energy), this energy is positive and one has the Cauchy-Schwarz inequality for the 
Dirichlet space :
\begin{displaymath} 
\bigg| \int_\D G_{\nu_1} \, d \nu_2 \bigg| \leq 
\bigg( \int_\D G_{\nu_1} \, d \nu_1 \bigg)^{1/2} \bigg(\int_\D G_{\nu_2} d \nu_2 \bigg)^{1/2} .
\end{displaymath} 
Applying this to the measures $\nu_1 = \mu$ and $\nu_2 = \nu = \Delta u. dz$, we get, since $G_\mu \leq 1$:
\begin{align*} 
\mu (K) 
& \leq \frac{1}{2 \pi} \bigg(\int_\D G_\mu (z) \, d \mu (z) \bigg)^{1/2} \bigg(\int_\D G_\nu (z) \Delta u (z) \, dz \bigg)^{1/2} \\ 
& \leq \frac{1}{2 \pi} \, [\mu (K)]^{1/2} \bigg(\int_\D G_\nu (z) \Delta u (z) \, dz \bigg)^{1/2} \\
& = \frac{1}{2 \pi} \, [\mu (K)]^{1/2} \bigg(\int_\D G_\nu \, d \nu \bigg)^{1/2}. 
\end{align*} 
Now, since $u \in {\cal C}_0^\infty (\D)$, one has G.~C.~Evans' theorem \cite{Evans} (see \cite{Brelot}, Chapitre~XI, Lemme~1, p.~141, or 
\cite{LAN}, Theorem~1.20, p.~97):
\begin{displaymath} 
\int_\D G_\nu \, d \nu  = 2 \pi \int_\D |\nabla u (z) |^2 \, dz \, .
\end{displaymath} 
Therefore, we get:
\begin{displaymath} 
\mu (K) \leq \frac{1}{2 \pi} \int_\D |\nabla u (z) |^2 \, dz \,. 
\end{displaymath} 
Taking the supremum on $\mu$ of the left-hand side and the infimum on $u$ of the right-hand side, we obtain:
\begin{displaymath} 
\capa (K) \leq \inf \Big\{\frac{1}{2 \pi} \int_\D |\nabla u (z) |^2 \, dz \, ; \  u \in {\cal C}_0^\infty (\D) \text{ and } u \geq 1 \text{ on } K \Big\} \, .
\end{displaymath} 
\par

2) Let $\eps > 0$. \par 
Let $K_j = \{ z \in \C \, ; \ {\rm dist}\, (z, K) \leq 1/j\}$, $j \geq 1$. Each $K_j$ is compact and is contained in $\D$ for $j$ large enough, say 
$j \geq j_0$. Since $K = \bigcap_{j \geq j_0} K_j$ (and the sequence is decreasing), one has 
$\capa (K_j) \mathop{\longrightarrow}\limits_{j \to \infty} \capa (K)$ 
(\cite{Conway}, Proposition~21.7.15; note that though this proposition is stated for the logarithmic capacity, the proof clearly works also for the Green capacity). 
Hence, there is some $j \geq j_0$ such that, for $K' = K_j$, one has $(1 + \eps) \, \capa (K) \geq \capa (K')$. \par

Let $\mu_0$ be an equilibrium measure of $K'$. One has $\mu_0 (K' ) = 1$, $I (\mu_0) = V (K')$, $G_{\mu_0} \leq V (K')$ on $\D$. Moreover, by 
\cite{Conway}, Lemma~21.10.1 (based on Frostman's theorem: see \cite{Conway}, Theorem~21.7.12, whose proof works also for the Green capacity), one has  
$G_{\mu_0} = V (K')$ on ${\rm int}\, (K')$, hence on $K$. Let $\mu = \capa (K') \, \mu_0$. Then $\mu (K') = \capa (K')$, 
$I (\mu) = [\capa (K')]^2 I (\mu_0) = \capa (K')$, and, since $G_\mu = \capa (K')\, G_{\mu_0}$, one has also 
$G_\mu \leq 1$ on $\D$ and $G_\mu = 1$ on $K$. \par 

By a theorem of G. Choquet \cite{Choquet2}, we can find, by regularization (\cite{Brelot}, p.~26 and Lemma, p.~135 and pp.~142--145, or \cite{LAN}, 
Theorem~1.9, p.~70, which applies since $G_\mu - U_2^\mu$ is a harmonic function) an increasing sequence of positive infinitely differentiable functions 
$v_n$ on $\D$ which converges pointwise to $G_\mu$ and such that:
\begin{displaymath} 
\int_\D |\nabla v_n (z)|^2 \, d z \mathop{\ \longrightarrow \ }_{n \to \infty} \int_ \D |\nabla G_\mu (z)|^2 \, d z \,. 
\end{displaymath} 
Since $(v_n)_n$ is increasing and converges pointwise to $1$ on the compact set $K$, Dini's theorem tells that one has uniform convergence. Hence, we can find 
some $v = v_n$ such that $v \geq (1 + \eps)^{- 1}$ on $K$ and 
\begin{displaymath} 
\int_\D |\nabla v (z)|^2 \, dz \leq (1 + \eps) \int_ \D |\nabla G_\mu (z) |^2 \, dz \,. 
\end{displaymath} 
Note that $v = 0$ on $\partial \D$ since $0 \leq v \leq G_\mu$, which is equal to $0$ on $\partial \D$. \par

Putting $u = (1 + \eps) v$, one has $u \in {\cal C}_0^\infty (\D)$, $u \geq 1$ on $K$ and 
\begin{displaymath} 
\int_\D |\nabla u (z)|^2 \, dz \leq (1 + \eps)^3 \int_ \D |\nabla G_\mu (z) |^2 \, dz \,. 
\end{displaymath} 

But we know by G.~C.~Evans's theorem (see \cite{Papa}, Proposition~7.3, or \cite{Brelot}, Chapitre~XI, p.~142 and pp.~144--145, or \cite{LAN}, 
Theorem~1.20, p.~97) that:
\begin{displaymath} 
I (\mu) = \frac{1}{2 \pi} \int_\D |\nabla G_\mu (z)|^2 \, dz \, .
\end{displaymath} 
We get hence:
\begin{align*}
(1 + \eps) \, \capa (K) 
& \geq \capa (K') = I (\mu) = \frac{1}{2 \pi} \int_\D |\nabla G_\mu (z) |^2 \, dz \\
& \geq \frac{1}{(1 + \eps)^3} \, \frac{1}{2 \pi} \int_\D |\nabla u (z) |^2 \, dz \,.
\end{align*}
Since $\eps > 0$ was arbitrary, we get:
\begin{displaymath} 
\capa (K) \geq \inf\Big\{\frac{1}{2 \pi} \int_\D |\nabla u (z) |^2 \, dz \, ; \  u \in {\cal C}_0^\infty (\D) \text{ and } u \geq 1 \text{ on } K \Big\} \, ,
\end{displaymath} 
and that ends the proof. \qed 
\goodbreak\bigskip

\noindent{\bf Remark.} After this paper was completed, we have found an alternative proof of Theorem~\ref{saksman}. We sketch it here. \par 
As in the above proof, we may assume that $0$ and $1 - \eps$ belong to $K$. Consider $K^\ast = \{ |z| \, ; \ z \in K\}$. Since $K$ is connected, the same 
holds for $K^\ast$. Hence the interval $[0, 1 - \eps]$ is contained in $K^\ast$. It follows that $\capa ([0, 1 - \eps]) \leq \capa (K^\ast)$. But we saw in 
Example~2 that $\capa ([0, 1 - \eps]) \approx \log (1 / \eps)$; hence $\capa (K^\ast) \gtrsim \log (1 / \eps)$. It remains to use that the map 
$\alpha \colon z \mapsto |z|$ is a contraction for the pseudo-hyperbolic metric and hence $\capa (K^\ast) \leq \capa (K)$ 
(see \cite{LAN}, Chap.~II, Theorem~2.9, and the comment p.~175 for the Green capacity). In fact, if $\nu$ is any probability measure supported by $K^\ast$, 
there exists (see \cite{harmonic}, Chap.~III, Lemma~4.6) a probability measure $\mu$ on $K$ such that $\alpha (\mu) = \nu$. Hence:
\begin{align*}
V (K) \leq I_K (\mu) 
& = \iint_{\D \times \D} g (z, w) \, d \mu (z) \, d \mu (w) = \iint_{\D \times \D} \log \frac{1}{\rho (z, w)} \, d \mu (z) \, d \mu (w) \\ 
& \leq \iint_{\D \times \D} \log \frac{1}{\rho (|z|, |w|)} \, d \mu (z) \, d \mu (w) \\
& = \iint_{\D \times \D} \log \frac{1}{\rho (z, w)} \, d \nu (z) \, d \nu (w) = I_{K^\ast} (\nu).
\end{align*}
Taking the infimum over all $\nu$, we get $V (K) \leq V (K^\ast)$. \qed 
\par
\bigskip

As a corollary of Theorem~\ref{saksman}, we get a new proof of \cite{LIQUEROD}, Theorem~3.4 and of \cite{LELIQURO}, Theorem~2.2. 

\begin{theorem} \label{coroll} 
There exists an absolute constant $c > 0$ such that, for any symbol $\varphi$ on a weighted analytic space $H$, one has:
\begin{displaymath} 
{\rm diam}_{\rho} \, [\varphi (\mathbb{D})] > r \quad \Longrightarrow \quad \beta (C_\varphi) \geq \exp \bigg[- \frac{c}{\log 1/(1 - r)} \bigg] \, .
\end{displaymath} 
In particular:
\begin{displaymath} 
\Vert \varphi \Vert_\infty = 1 \quad \Longrightarrow \quad \beta (C_\varphi) = 1 \, .
\end{displaymath} 
\end{theorem}

\noindent{\bf Proof.} The first statement is a direct consequence of Theorem~\ref{principal}, modulo Theorem~\ref{alano} and Theorem~\ref{saksman},  
applied to $\varphi (\mathbb{D})$ and its closure. \par 
One cannot replace ${\rm diam}_{\rho} \, [\phi (\D)] > r$ by $\Vert \varphi \Vert_\infty > r$ in this first statement as indicated by the following example: 
\begin{displaymath} 
\phi (z) = \frac{a - (z/2)}{1 - \overline{a} (z/2)} = \Phi_a [h (z)] \, , 
\end{displaymath} 
where $\Phi_{a} (z) = \frac{a - z}{1 - \overline{a} z}$ with $a \in \mathbb{D}$ and $h (z) = z/2$ is the dilation with ratio $1/2$. Then 
$\Vert \varphi \Vert_\infty \geq |\Phi_{a}(0)|=|a|$ and $\beta (C_\varphi) = \beta (C_h) = 1/2$. \par

However, one can do so if moreover $\varphi (0) = 0$ because then, clearly: 
\begin{displaymath} 
\Vert \varphi \Vert_\infty > r \quad \Longrightarrow \quad {\rm diam}_{\rho} \, [\varphi (\mathbb{D})] > r \, .
\end{displaymath} 
This is enough for the second statement since, putting $a = \varphi (0)$, we have, due to the fact that $\Phi_a$ is unimodular on the 
whole unit circle: $\Vert \Phi_a \circ \varphi \Vert_\infty = \Vert \varphi \Vert_\infty = 1$, 
$(\Phi_a \circ \varphi) (0) = 0$ and $\beta (C_\varphi) = \beta (C_{ \Phi_a \circ \varphi})$. \qed

\subsection{A remark} 

We proved in \cite{Dirichlet} that every composition operator $C_\phi$ which is bounded on the Dirichlet space ${\cal D}$ is compact on the Hardy space 
$H^2$ (and hence on the Bergman space ${\mathfrak B}^2$), and even in all Schatten classes on $H^2$ and ${\mathfrak B}^2$. So one may expect that the 
approximation numbers of composition operators on the Dirichlet space are bigger than those on the Hardy space (and bigger than those on the Bergman space). 
Since Theorem~\ref{principal} and Theorem~\ref{ancien} show that $\beta (C_\phi)$ is the same for these three spaces, it follows that the answer will be 
certainly quite subtle and cannot only involve $\log a_n (C_\phi)$.

\section{The $H^p$ case, $1 \leq p < \infty$} 

Here, we consider the case of composition operators on $H^p$ for $1 \leq p < \infty$. 
\par\medskip

For every $a \in \D$, we denote by $e_a \in (H^p)^\ast$ the evaluation map at $a$, namely:
\begin{equation} 
\qquad e_a (f) = f (a) \, , \quad f \in H^p .
\end{equation} 
We know that (\cite{ZHU}, p.~253):
\begin{equation}\label{eval} 
\| e_a \| = \left(\frac{1}{1 - |a|^2} \right)^{1/p}  
\end{equation}
and the mapping equation 
\begin{equation} 
C_{\varphi}^{\ast} (e_a) = e_{\varphi(a)} 
\end{equation} 
still holds. \par 
\medskip

Throughout this section we denote by $\| \, . \, \|$, without any subscript, the norm in the dual space $(H^p)^\ast$. \par
\smallskip 

Let us stress that this dual norm of $(H^p)^\ast$ is, for $1 < p < \infty$, equivalent, but not equal, to the norm $\| \, . \, \|_q$ of $H^q$, and the 
equivalence constant tends to infinity when $p$ goes to $1$ or to $\infty$. 

With this preliminaries, we are going to see that Theorem~\ref{principal} remains true. \par
\smallskip

\begin{theorem} \label{encore}
Let $1 \leq p < \infty$ and $C_\phi \colon H^p \to H^p$. \par

$1)$ If $\overline{\phi (\D)} \subseteq \D$, then:
\begin{displaymath} 
\beta (C_\phi) = \e^{- 1 / \capa [\phi (\D)]} \, .
\end{displaymath} 
\par

$2)$ One has: 
\begin{displaymath} 
\| \phi \|_\infty = 1 \quad \Longrightarrow \quad \beta (C_\phi) = 1 \, .
\end{displaymath} 
\end{theorem}
\goodbreak

We begin with the following lemma, which extends Lemma~\ref{known}.  

\begin{lemma} \label{caspe}  
Let $X$ be a Banach space, and  $T \colon X \to X$ be a compact operator. Let us set:
\begin{equation} \label{andrei} 
\varepsilon_{n} (T) = \inf_{\dim E < n} \left[\sup_{x \in B_X} {\rm dist}\, (T x, T E) \right]. 
\end{equation}
Then $\varepsilon_{n} (T)\leq 2 \, \sqrt{n} \, c_{n} (T)$.
\end{lemma}

\noindent {\bf Proof.} Let $\varepsilon > 0$, and let $F$ be a subspace of $X$ of codimension $ < n$ such that 
$\Vert T_{\mid F} \Vert \leq c_{n} (T) + \varepsilon$. Let $Q \colon X \to F$ be an onto projection of norm  
$\Vert Q \Vert \leq \sqrt{n}+ 1 \leq 2 \, \sqrt n$ (see \cite{LQ}, Chapitre 5, Th\'eor\`eme~III. 4, 2), or \cite{Wojtaszczyk}, III.B.11) and let 
$R = T (I - Q)$. Then $E = (I - Q) X$ satisfies $\dim E < n$. If $x \in B_X$, the closed unit ball of $X$, then:
\begin{displaymath} 
{\rm dist}\, (T x, T E) \leq \Vert T x - R x \Vert = \Vert T Q x \Vert \leq \Vert T_{\mid F} \Vert \, \Vert Q x\Vert 
\leq  (c_{n} (T) + \varepsilon) \, 2 \, \sqrt n \, .
\end{displaymath} 
This implies $\varepsilon_{n} (T) \leq 2 \, \sqrt n \, (c_{n}(T) + \varepsilon)$. \par \smallskip

The result follows since $\varepsilon$ was arbitrary. \qed
\bigskip

\noindent {\bf Proof of Theorem~\ref{encore}.} $1)$ a) We first prove that $\beta^- (C_\phi) \geq \e^{- 1 / \capa [\phi (\D)]}$. \par 
Let $\tilde L_r = \sup_{|a| \leq r} \| e_a \| = \big(\frac{1}{1 - r^2} \big)^{1/p}$, for $0 < r < 1$. Using the same notations and estimations as in 
Theorem~\ref{below}, up to the replacement of $L_r$ by $\tilde L_r$, we get:
\begin{displaymath} 
\varepsilon_{n} (T) \geq (1 - \varepsilon) \, \tilde L_{r_j}^{- 1} \, \alpha \,  [M (K_j)]^n \, . 
\end{displaymath} 
Lemma~\ref{caspe} now implies: 
\begin{displaymath} 
a_{n} (T) \geq c_{n} (T) \geq \alpha \, \frac{1 - \varepsilon}{2 \, \sqrt n} \, \tilde L_{r_j}^{- 1} [M (K_j)]^n \, . 
\end{displaymath} 
The rest of the proof is unchanged, since the presence of the factor $1/\sqrt n$  does not affect the result. \par\smallskip

\quad b) The upper bound is even simpler since $H^\infty \subseteq H^p$. For example, with the notations of Section~\ref{section upper}, setting 
$A (f) = h \circ \varphi$ as in \eqref{choice}, we can replace Lemma~\ref{ext} by 
\begin{displaymath} 
\Vert g \circ \varphi - h \circ \varphi \Vert_p \leq \Vert g \circ \varphi - h \circ \varphi\Vert_\infty = \Vert g - h \Vert_{\mathcal{C}(K)} \, ,
\end{displaymath} 
where $K = \overline{\varphi(\D)}$. \par \smallskip

2) That follows from Theorem~\ref{saksman}, as in Section~\ref{saksman}. 
\qed

\bigskip

\noindent{\bf Ackowledgements.} We thanks A. Ancona and E. Saksman for crucial informations on the Green capacity. \par\bigskip
\goodbreak


\bigskip

\noindent
{\rm Daniel Li}, Univ Lille Nord de France, \\
U-Artois, Laboratoire de Math\'ematiques de Lens EA~2462 \\ 
\& F\'ed\'eration CNRS Nord-Pas-de-Calais FR~2956, \\
Facult\'e des Sciences Jean Perrin, Rue Jean Souvraz, S.P.\kern 1mm 18, \\
F-62\kern 1mm 300 LENS, FRANCE \\ 
daniel.li@euler.univ-artois.fr
\medskip

\noindent
{\rm Herv\'e Queff\'elec}, Univ Lille Nord de France, \\
USTL, Laboratoire Paul Painlev\'e U.M.R. CNRS 8524 \& 
F\'ed\'eration CNRS Nord-Pas-de-Calais FR~2956, \\
F-59\kern 1mm 655 VILLENEUVE D'ASCQ Cedex, 
FRANCE \\ 
Herve.Queffelec@univ-lille1.fr
\smallskip

\noindent
{\rm Luis Rodr{\'\i}guez-Piazza}, Universidad de Sevilla, \\
Facultad de Matem\'aticas, Departamento de An\'alisis Matem\'atico \& IMUS,\\ 
Apartado de Correos 1160,\\
41\kern 1mm 080 SEVILLA, SPAIN \\ 
piazza@us.es

\end{document}